\documentclass[11pt,amsfonts]{article}
\usepackage{graphicx}
\usepackage{latexsym}
\usepackage{amssymb}
\usepackage{amsmath}
\usepackage{enumerate}
\usepackage{color}
\usepackage{layout}
\usepackage{eufrak}

\newtheorem{prop}{Proposition}
\newtheorem{lemma}{Lemma}

\newtheorem{corollary}{Corollary}
\newtheorem{theorem}{Theorem}
\newtheorem{remark}{Remark}

\def\real{{\mathord{{\rm I\kern-2.8pt R}}}}        
\def\inte{{\mathord{{\rm I\kern-2.8pt N}}}}

\def\Dom{{\mathrm{{\rm Dom}}}}

\def\sZZ{{\rm Z\kern-2.8ptem{}Z}}

\def\z{{\mathchoice
  {\sZZ}
  {\sZZ}
  {\rm Z\kern-0.30em{}Z}
  {\rm Z\kern-0.25em{}Z} }}
\def\sQQ{{\kern 0.27em \vrule height1.45ex width0.03em depth0em
          \kern-0.30em \rm Q}}
\def\qu{{\mathchoice
    {\sQQ}
    {\sQQ}
  {\kern 0.225em \vrule height1.05ex width0.025em depth0em \kern-0.25em \rm Q}
  {\kern 0.180em \vrule height0.78ex width0.020em depth0em \kern-0.20em \rm Q}
        }}
\def\sCC{{\kern 0.27em \vrule height1.45ex width0.03em depth0em
          \kern-0.30em \rm C}}
\def\complex{{\mathchoice
    {\sCC}
    {\sCC}
  {\kern 0.225em \vrule height1.05ex width0.025em depth0em \kern-0.25em \rm C}
  {\kern 0.180em \vrule height0.78ex width0.020em depth0em \kern-0.20em \rm C}
        }}

\newcommand{\ba}{\begin{array}}
\newcommand{\ea}{\end{array}}
\newcommand{\be}{\begin{equation}}
\newcommand{\ee}{\end{equation}}
\newcommand{\bea}{\begin{eqnarray}}
\newcommand{\eea}{\end{eqnarray}}
\newcommand{\beaa}{\begin{eqnarray*}}
\newcommand{\eeaa}{\end{eqnarray*}}

\def\b{\beta}

\def\z{\zeta}

\font\tenmath=msbm10 \font\sevenmath=msbm7 \font\fivemath=msbm5
\newfam\mathfam \textfont\mathfam=\tenmath
\scriptfont\mathfam=\sevenmath \scriptscriptfont\mathfam=\fivemath

\def \b{\noindent}

\def \={{\buildrel {\rm (law)} \over =}}

\def\qed{ \hfill \vrule width.25cm height.25cm depth0cm\smallskip}

\newcommand{\basa}{\begin{assumption}}
\newcommand{\easa}{\end{assumption}}

\newcommand{\bas}{\begin{assum}}
\newcommand{\eas}{\end{assum}}

\newcommand{\ignore}[1]{}
\begin{document}

\renewcommand{\thefootnote}{\fnsymbol{footnote}}

\renewcommand{\thefootnote}{\fnsymbol{footnote}}

\title{Parameter estimation for the Rosenblatt Ornstein-Uhlenbeck process with periodic mean }
\author{R. Shevchenko$^{\dagger}$  and C. A. Tudor$^{*}$\vspace*{0.2in} \\
 $^{\dagger}$ Fakult\"at f\"ur Mathematik,
LSIV, TU Dortmund \\
 Vogelpothsweg 87, 44227 Dortmund, Germany. \\
\quad radomyra.shevchenko@tu-dortmund.de \vspace*{0.1in}\\
$^{*}$ Laboratoire Paul Painlev\'e, Universit\'e de Lille \\
 F-59655 Villeneuve d'Ascq, France.\\
 \quad ciprian.tudor@math.univ-lille.fr\\
\vspace*{0.1in}}

\maketitle

\begin{abstract}
We study the least squares estimator for the drift parameter of the Langevin stochastic equation driven by the Rosenblatt process. Using the techniques of the Malliavin calculus  and the stochastic integration with respect to the Rosenblatt process, we analyze the consistency and the asymptotic distribution of this estimator. We also introduce alternative estimators, which can be simulated, and we study their asymptotic properties.
\end{abstract}

\vskip0.3cm

{\bf 2010 AMS Classification Numbers:}  60H15, 60H07, 60G35.

 \vskip0.3cm


\section{Introduction}

While the parameter estimation for continuously observed classical diffussion processes has a long history (see e.g. \cite{Ku} and the references therein), the statistical inference for stochastic equations driven by fractional Brownian motion (fBm) and related processes started more recently, in the nineties. Since then a large number of reserarch articles considered the problem of drift parameter estimation for various fractional diffussions and in particular for the fractional Ornstein-Uhlenbeck process, which is defined as the solution to the Langevin equation
\begin{equation}
\label{intro-1}
dX_{t}= -\alpha X_{t} dt+ \sigma dB ^{H}_{t}, \hskip0.5cm t\geq 0,
\end{equation}
with $\alpha \in \mathbb{R}, \sigma >0$ and $( B ^{H}_{t}) _{t\geq 0} $ a fractional Brownian motion with Hurst parameter $H\in (0, 1). $ We refer among many others to 
 \cite{KL}, \cite{Rao1},  \cite{HN} or \cite{TV}. More recently, a non-Gaussian extension of the model (\ref{intro-1})  has been considered by several authors (see e.g. \cite{ntra}, \cite{SlaTud}), by replacing the fractional Brownian noise in (\ref{intro-1}) by a Hermite process.  The Hermite processes are self-similar processes with stationary increments and long memory, with the same covariance as the fBm, but non-Gaussian. The Hermite process of order $q\geq 1$ lives in the $q$th Wiener chaos, i.e. it can be expressed as a iterated stochastic integral with respect to the Wiener process. For $q=1$ it coincides with the fBm (which is the only Gaussian Hermite process) while for $q=2$ it is known as the Rosenblatt process. 

We will consider the following model: 
\begin{equation}
\label{intro-2} 
dX_{t} = (L (t)- \alpha X _{t}) dt + d Z ^{H} _{t}, \hskip0.5cm t\geq 0,
\end{equation}
where the random noise $(Z ^{H}_{t}) _{t\geq 0} $ is a Rosenblatt process with self-similarity order $ H\in \left( \frac{1}{2}, 1\right) $ and $L$ is a periodic function. We will assume that $L$ can be written as $L(t) = \sum_{i=1}^{p} \mu _{i} \varphi_{i} (t)$ with some suitable  known  periodic functions $\varphi$, $i=1,..,p$. The purpose is to estimate jointly the parameters $\mu_{1},.., \mu_{p}$ and $\alpha $ based on a continuous-time observation of the solution to (\ref{intro-2}). Models similar to (\ref{intro-2}) have been considered in \cite{FK} for the case of Wiener noise and in \cite{BESV}, \cite{DFW} for fractional Brownian noise. These models are proposed in order to better capture the characteristics of the empirical data in some applications (related to seasonalities, for example).

We estimate the parameters  $\mu_{1},.., \mu_{p}$ and $\alpha $ in (\ref{intro-2})  by using a least-square estimator introduced in \cite{FK} or \cite{DFW}.  The resulting estimator involves stochastic integrals with non-deterministic integrands with respect to the Rosenblatt process and this fact makes  its analysis more complex. Since the Rosenblatt process is neither a semimartingale nor a Gaussian process, we cannot use the classical stochastic integration with respect to it.  Instead, we will use the stochastic analysis of the Rosenblatt process developed in \cite{T}. We show the consistency of the estimator and we find its asymptotic behavior in distribution. Our proofs are based on the Malliavin calculus, the correlation structure of the solution to (\ref{intro-2}) and the properties of the random variables living in the second Wiener chaos.

We organized our paper as follows. In Section 2 we introduce some preliminaries concerning the Rosenblatt process and in Section 3 we discuss the details of our setting and demonstrate some auxiliary results. Sections 4 and 5 are concerned with the construction of the least squares estimator as well as with the analysis of its asymptotic behaviour. In Section 6, alternative estimators are defined and studied. Finally, the Appendix serves as a brief compendium of definitions and results from Malliavin calculus.

\section{Preliminaries: The Rosenblatt process and the stochastic integral with respect to it}\label{2.1}

Let us start by recalling the definition and the basic properties of the Rosenblatt process as well as the construction of the stochastic integral with respect to this process, which is neither Gaussian nor a semimartingale. For a more complete exposition, we refer to the monographs \cite{PiTa-book},
 \cite{T1} or to the reference \cite{T}.  Notice that there are  several possibles definitions of the Rosenblatt process. Here, we chose to work with the so-called finite interval representation of it. Let $H>\frac{1}{2}$ and $(B_{t})_{t\geq 0} $ a Brownian motion. Consider the kernel
\begin{equation}
\label{K}
K ^{H}(t,s)= c_{H} s ^{\frac{1}{2}-H} \int_{s} ^{t} (u-s) ^{H-\frac{3}{2}}u ^{H-\frac{1}{2}} du
\end{equation}
with $t>s$ and $c_{H}$ a deterministic constant.{\it  The Rosenblatt process with self-similarity index $H\in \left(\frac{1}{2}, 1\right)$ }  is defined as 
\begin{equation}
\label{rose}
Z ^{H}_{t}= d(H) \int_{0} ^{t} \int_{0} ^{t} \left( \int_{y_{1}\vee y_{2}}^{t} \frac{\partial K ^{ H'}}{\partial u} (u, y_{1})\frac{\partial K ^{ H'}}{\partial u} (u, y_{2})du\right) dB_{y_{1}} dB_{y_{2}}, \hskip0.3cm t\geq 0
\end{equation}
with $$H'=\frac{H+1}{2}$$ and $d(H)$ a deterministic constant that ensures $\mathbf{E}  (Z_{t} ^{H}) ^{2}= t^{2H}$ for every $t\geq 0$.  The stochastic integral in (\ref{rose}) is a multiple integral of order 2 with respect to the Wiener process $B$, see the Appendix. The process $\left( Z ^{H}_{t}\right)_{t\geq 0}$ is a self-similar stochastic process (with the self-similarity index $H$) with stationary increments, living in the second Wiener chaos, with H\"older continuous paths of order $\delta \in (0, H)$. 

Let us denote by $\mathcal{H}$ the canonical Hilbert space associated to the fractional Brownian motion with parameter $H$, i.e. $\mathcal{H}$ is the closure of the linear space generated by the indicator functions $\{ 1_{[0, t] }, t\geq 0\} $ with respect to the inner product
\begin{equation*}
\langle 1_{[0,t]}, 1_{[0, s]} \rangle _{\mathcal{H}}=\frac{1}{2} \left( t^{2H}+s^{2H} -\vert t-s\vert ^{2H} \right), \hskip0.3cm t,s \geq 0. 
\end{equation*}
It is also possible to define Skorohod integrals of random integrands with respect to the Rosenblatt process. For a square integrable stochastic process $(g_{t})_{t\geq 0}$ we set 
\begin{equation}
\label{30i-3}
\int_{0} ^{T} g_{s} dZ^{H} _{s}:= \int_{0} ^{t} \int_{0} ^{t} I(g) (y_{1}, y_{2})dB_{y_{1}} dB_{y_{2}} 
\end{equation}
with  the transfer operator
\begin{equation}
\label{I}
I(g)(y_{1}, y_{2}) = \int_{ y_{1} \vee y_{2}} ^{T} g_{u} \frac{ \partial K ^{ H'} }{\partial u} (u, y_{1})\frac{ \partial K ^{ H'} }{\partial u} (u, y_{2}) du.
\end{equation}
The notation  $dB$  in (\ref{30i-3})  indicates the Skorohod integral with respect to the Wiener process $(B_{y})_{y\geq 0}$. From Lemma 1 in \cite{T}, the Skorohod integral (\ref{30i-1}) is well-defined if 
\begin{equation}\label{30i-4}
\mathbf{E} \int_{0} ^{T} \int_{0} ^{T} \Vert D_{x_{1}, x_{2}} g\Vert ^{2}_{\mathcal{H}}dx_{1} dx_{2} <\infty.
\end{equation}
Moreover, if $g\in \mathbb{L } ^{2, p}$ ($p\geq 2$), then for every $t\geq 0$

\begin{equation}
\label{30i-5}
\mathbf{E}  \left| \int_{0} ^{t} g_{s} d Z^{H}_{s} \right| ^{p} \leq c(p, H) \sup_{ t\in [0, t]} \left[ \mathbf{E}  \vert g_{r} \vert ^{p} + \mathbf{E}  \int_{0} ^{t} \int_{0} ^{t} \Vert D ^{(2)}_{x_{1}, x_{2}} g_{r} \Vert ^{p} dx_{1}dx_{2} \right] t ^{pH}. 
\end{equation}
If $g\in \mathcal{H}$ is deterministic, then the integral (\ref{30i-3}) is a Wiener integral  with respect to the Rosenblatt process (also called Wiener-Rosenblatt integral) and it satisfies the following isometry
\begin{equation*}\label{isoZ}
\mathbf{E}  \left( \int_{0} ^{t} g_{s} d Z^{H}_{s} \int_{0} ^{t} h_{s} d Z^{H}_{s}\right)=H(2H-1) \int_{0} ^{t} \int_{0} ^{t} g(u) h(v) \vert u-v\vert ^{2H-2}  =\:langle g, h\rangle _{\mathcal{H}}
\end{equation*}
for ant $t\geq 0$ and for any functions $g,h$ such that $\int_{0} ^{t} \int_{0} ^{t} \vert g(u) h(v)\vert  \vert u-v\vert ^{2H-2} du dv  <\infty.$

\section{The Rosenblatt Ornstein-Uhlenbeck process with periodic mean}

The Rosenblatt Ornstein-Uhlenbeck (ROU in the sequel) process is defined as the solution of the Langevin equation driven by  a Rosenblatt noise, see e.g. \cite{tma} or \cite{SlaTud}.  The ROU process with periodic mean is defined as the solution to the Langevin equation whose drift is a periodic function. More precisely, we will consider the stochastic differential equation 
\begin{equation}
\label{sde}
X_t=\int_0^t \left( L(s)-\alpha X_s \right)ds + Z^H_t, \hskip0.5cm t\geq 0
\end{equation}
with vanishing initial condition, where  $Z^H$ is the Rosenblatt process with self-similarity index $H\in \left(\frac{1}{2}, 1\right)$. $L$ is assumed to be a deterministic function that can be expressed as a linear combination of known bounded $1$-periodic functions (assumed to be orthonormal in $L^ {2} ([0,1])$,  without loss of generality), i.e., for $p\geq 1$,
\begin{equation}
\label{L}
L(s)=\sum_{i=1}^p \mu_i\varphi_i(s),\quad s\geq 0.
\end{equation}

Let us focus on the basic properties of the solution to (\ref{sde}). As in the case when the noise is a fractional Brownian motion, it can be shown that (\ref{sde}) admits a unique strong solution which can be written as 
\begin{equation}\label{x}
X_t=e^{-\alpha t}\left(\int_0^t e^{\alpha s}L(s) ds + \int_0^t e^{\alpha s}dZ^H_s\right)=:h(t)+Y_t, \hskip0.4cm t\geq 0,
\end{equation}
where we use the notation
\begin{equation}
\label{hy}
h(t)= e^{-\alpha t} \int_0^t e^{\alpha s}L(s) ds \mbox{ and } Y_{t}= e^{-\alpha t} \int_0^t e^{\alpha s}dZ^H_s 
\end{equation}
for every $t\geq 0$. The stochastic integral $dZ ^{H}_{s} $ in (\ref{x}) is a Wiener integral with respect to the Rosenblatt process $ Z^{H}$ and  we will call the process $(X_{t})_{t\geq 0} $ {\it the Rosenblatt Ornstein Uhlebeck process with periodic mean}.  We can also define the so-called stationary Rosenblatt Ornstein-Uhlenbeck process with periodic mean by putting 

\begin{equation}
\label{tx}
\tilde{X}_t=e^{-\alpha t}\left(\int_{-\infty}^t e^{\alpha s}L(s) ds + \int_{-\infty}^t e^{\alpha s}dZ^H_s\right)=:\tilde{h}(t)+\tilde{Y}_t, \hskip0.4cm t\geq 0,
\end{equation}
with
\begin{equation}
\label{thy}
\tilde{h}(t)= e^{-\alpha t}\int_{-\infty}^t e^{\alpha s}L(s) ds \mbox{ and } \tilde{Y}_{t}=e^{-\alpha t}\int_{-\infty}^t e^{\alpha s}dZ^H_s.
\end{equation}
The existence of the stochastic integrals in (\ref{x}) and (\ref{tx}) has been showed in e.g. \cite{CKM} or \cite{tma}. We also recall the correlation structure of the process $\tilde{Y}$ (see \cite{CKM} or \cite{tma}): for every $t\geq 0$ and for $s\to \infty$ we have with $c_{H}\in \mathbb{R}$ 
\begin{equation}
\label{30i-1}
\mathbf{E}  \tilde{Y} _{t} \tilde{Y}_{t+s}= c_{H} s ^{2H-2}+ O (s ^{2H-4})
\end{equation}

We will start by proving some ergodic type properties of the process $X$. These properties will be needed in order  to analyze the asymptotic properties of our estimators in the sequel.

\begin{prop}\label{pp1}
Let $\varphi: \mathbb{R} \to \mathbb{R}$ be a bounded  1-periodic function and let $(\tilde{Y}_{t}) _{t\geq 0} $ be given by (\ref{thy}). Then 
 \[\frac{1}{n} \int_0^n \varphi (t) \tilde{Y}_t dt\to _{n \to \infty}  0 \mbox{ almost surely. }\]

\end{prop}
{\bf Proof: } We have for every $n\geq 1$
$$ \mathbf{E} \left[\left(\frac{1}{n} \int_0^n \varphi (t) \tilde{Y}_t dt \right)^2\right]=\frac{1}{n^2}\int_0^n\int_0^n \varphi (t)\varphi (s) \mathbf{E}  [\tilde Y_t \tilde Y_s]dtds.
$$
First notice that for every integer $n_{0}< n$ we have 
\begin{equation}
\label{30i-2}
 \frac{1}{ n ^{2}} \int_{ [0, n] ^{2}\setminus [ n_{0}, n ] ^{2}}  \varphi (t)\varphi (s) \mathbf{E}  [\tilde Y_t \tilde Y_s]dtds \to _{n\to \infty }0.
\end{equation}
Indeed, using the notation $a_{n} \lesssim b_{n}$ to indicate that for $n$ large we have $a_{n} \leq c b_{n} +c_{n} $ with $c_{n} \to _{n\to \infty} 0$, we can write
\begin{eqnarray*}
 &&\frac{1}{ n ^{2}} \int_{ [0, n] ^{2}\setminus [n_{0}, n ] ^{2}}  \varphi (t)\varphi (s) \mathbf{E}  [\tilde Y_t \tilde Y_s]dtds \\
&=& \frac{1}{n ^{2}} \int_{0} ^{n_{0} } \int_{0} ^{n_{0}}  \varphi (t)\varphi (s) \mathbf{E}  [\tilde Y_t \tilde Y_s]dtds +2\frac{1}{n ^ {2}} \int_{0} ^{n_{0}} \int_{ n_{0}} ^{n}  \varphi (t)\varphi (s) \mathbf{E}  [\tilde Y_t \tilde Y_s]dtds\\
&\lesssim & \frac{1}{n ^{2}} \int_{0} ^{n_{0}} \int_{ n_{0}} ^{n}  \vert \varphi (t)\varphi (s) \vert ( \tilde{Y}_{t} ^{2}+ \tilde{Y}_{s} ^{2})dt ds\leq c n ^{-1},
\end{eqnarray*}
where we used  $\mathbf{E}  Y_{t} ^{2} \leq c $ for every $t\geq 0$ (see relation (2.16) in \cite{ntra}). We obtain by (\ref{30i-2}) and the periodicity of $\varphi$
\begin{eqnarray*}
\mathbf{E}  \left[\left(\frac{1}{n} \int_0^n \varphi (t) \tilde{Y}_t d \right)^2\right]
&\lesssim& \frac{1}{n^2}\int_{n_0}^n\int_{n_0}^n |\varphi (t)\varphi (s)| |t-s|^{2H-2}dtds\\
&\lesssim& \frac{1}{n^2}\int_{0}^n\int_{0}^n |\varphi (t)\varphi (s)| |t-s|^{2H-2}dtds\\
&\lesssim &  \frac{1}{n^2}\sum_{i,\,j = 0}^{n-1}\int_0^1\int_0^1 |\varphi (t)\varphi (s)| |t-s- (i-j)|^{2H-2}dtds\\
&\lesssim & \frac{1}{n^2}\sum_{i,\,j = 0; \vert i-j\vert < 2}^{n-1}\int_0^1\int_0^1 |\varphi (t)\varphi (s)| | (i-j)-(t-s)|^{2H-2}dtds  \\
&&\qquad +2 \frac{1}{n^2}\sum_{i,\,j = 0; i-j \geq 2}^{n-1}\int_0^1\int_0^1 |\varphi (t)\varphi (s)| ( (i-j)-(t-s))^{2H-2}dtds .
\end{eqnarray*}
We obtain
\begin{eqnarray*}
&&\frac{1}{n^2}\sum_{i,\,j = 0; \vert i-j\vert < 2}^{n-1}\int_0^1\int_0^1 |\varphi (t)\varphi (s)| | (i-j)-(t-s)|^{2H-2}dtds\\
&\lesssim & \frac{1}{n^2}n \max \left(\int_0^1\int_0^1 |\varphi (t)\varphi (s)| |t-s|^{2H-2}dtds,\, \int_0^1\int_0^1 |\varphi (t)\varphi (s)| | 1-(t-s)|^{2H-2}dtds\right).
\end{eqnarray*}
Because $ \varphi $ is bounded and $H>\frac{1}{2}$,  the two integrals above are finite and then the summand converges to zero as $n\to \infty$.

For the second summand note that
\[ ( (i-j)-(t-s))^{2H-2} =\left(1-\frac{t-s}{i-j}\right)^{2H-2} (i-j)^{2H-2},\]
and since for $ i-j \geq 2$ we have $ 1-\frac{t-s}{i-j} \geq \frac{1}{2}$, we deduce that this summand is bounded by
\[\frac{1}{n^2}\sum_{i,\,j = 0; \vert i-j\vert \geq 2}^{n-1}\int_0^1\int_0^1 |\varphi (t)\varphi (s)| |i-j|^{2H-2}dtds\]
up to a constant. In total, we have
\begin{eqnarray*}
\mathbf{E}  \left[\left(\frac{1}{n} \int_0^n \varphi (t) \tilde{Y}_t d \right)^2\right]
&\lesssim&\frac{1}{n^2}\sum_{i,\,j = 0; \vert i-j\vert \geq 2}^{n-1}\int_0^1\int_0^1 |\varphi (t)\varphi (s)| |i-j|^{2H-2}dtds\\
&\lesssim & \|\varphi\|_{L^2([0,\,1])}^2 \frac{1}{n^2} \sum_{i,\,j = 0}^{n-1} |i-j|^{2H-2}\lesssim n^{2H-2}.\\
\end{eqnarray*}
Since $\tilde Y_t$ is a second Wiener chaos element, so is the integral $ \int_0^n \varphi (t) \tilde{Y}_t dt$ as a pointwise limit. Therefore, due to  the hypercontractivity  property (\ref{hypercontractivity}) we obtain the bound
\[\mathbb E \left[\left(\frac{1}{n} \int_0^n \varphi (t) \tilde{Y}_t d \right)^{2m}\right]\lesssim n^{m(2H-2)}.\]
 We can choose an $m\in\mathbb N$  big enough, depending on $H$, such that the statement follows by the usual Borel-Cantelli argument.
\qed

 As a consequence of Proposition \ref{pp1},  we can deduce a discrete ergodic property for the shifted process $\tilde{X}$.

\begin{corollary}\label{r1}
For every $n\geq 1$, define  the process   ${\bf Y}_n:=\{\tilde{Y}_{n+s},\, s\in [0,\,1]\}$. Then ${\bf Y }$  satisfies the following discrete ergodic property 
\[ \frac{1}{n} \sum_{i=0}^{n-1}\int_0^{1}\varphi (t){\bf Y}_i(t)dt\to _{n\to \infty}0 \mbox{ almost surely.}
\]

Moreover, the process ${\bf X}_n:=\{\tilde{X}_{n+s},\, s\in [0,\,1]\}$ ($n\in \mathbb N$) also satisfies the discrete ergodic property, i.e.
\begin{equation*}
\frac{1}{n} \sum_{i=0}^{n-1}\int_{0}^{1}\varphi (t){\bf X}_i(t)dt\to _{n\to \infty}\int_{0} ^{1} \varphi(t) \tilde{h} (t)dt \mbox{ almost surely.}
\end{equation*}
\end{corollary}
{\bf Proof: } For ${\bf Y}_{n}$, the conclusion follows since
 $$\frac{1}{n} \sum_{i=0}^{n-1}\int_0^{1}\varphi (t){\bf Y}_i(t)dt= \frac{1}{n} \sum_{i=0} ^{n-1} \int_{i}^{i+1} \varphi(t) \tilde{Y}_{t}dt = \frac{1}{n} \int_{0} ^{n} \varphi (t) \tilde{Y}_{t}dt $$
while  for ${\bf X} _{n}$  we simply use the fact that $\tilde{h}$ is $1$-periodic. \qed

\section{The least squares estimator}
We will analyze the least squares estimator for the parameters of the model (\ref{sde}), inspired by the construction in \cite {FK} and \cite{DFW}. In the first part we recall its definition and basic  properties and in the second part we study its consistency and its limit behavior in distribution.

\subsection{Definition and basic properties}
Our purpose is to estimate the $(p+1)$-dimensional parameter 
\begin{equation}
\label{18f-2}
\vartheta =\left( \mu _{1},..., \mu _{p}, \alpha\right) 
\end{equation}
where $\mu _{i}$, $i=1,.., p$, are the coefficients that appear in the definition of the periodic function $L$ (see formula (\ref{L})) while $\alpha $ is the drift parameter of the ROU process  (\ref{sde}).  We will construct a least squares estimator (LSE) to estimate $\vartheta$. The construction of this estimator borrows the idea from 
\cite{FK} and \cite{DFW}. In these references, a general definition of the LSE of the parameter $\theta \in \mathbb{R} ^{p+1}$ of the stochastic differential equation
\begin{equation*}
dX_{t} =\theta f(t, X_{t}) dt + dB_{t}, \hskip0.5cm t \in [0, T] 
\end{equation*}
is presented,  $B$ being a general  noise. The idea is to minimize the error function 
$$\theta \to \sum_{i=1} ^{N} \left( X_{(i+1) \Delta t} -X_{i\Delta t} -\sum_{j=1}^{p+1}   f_{j} (i\Delta t, X_{i\Delta t}) \theta _{j} \Delta t \right) ^{2} $$
where $t_{i}, i=1,.., N$, denotes an equidistant  discretization of $[0, T] $ with $\Delta t= \frac{T}{N} $ and $f_{j}$ are the components of the function $f: \mathbb{R}_{+}\times \mathbb{R} ^{p+1}\to \mathbb{R}. $ As in \cite{FK}, \cite{DFW} we are led to the following LSE 
\begin{equation}
\label{tn}
\hat{\vartheta}_n:=(\hat{\mu}^1_n,\dots, \hat{\mu}^{p}_n,\,\hat \alpha_n):=Q_n^{-1}P_n,
\end{equation}
with  the $(p+1)$- dimensional random vector $P_{n}$  given by  ("{T}" denotes the transpose)
\begin{equation}
\label{pn}
P_n:=\left(\int_0^{n} \varphi_1 (t) dX_t,\dots ,\int_0^{n} \varphi_p (t) dX_t, -\int_0^{n} X_t dX_t\right)^T
\end{equation}
and with the matrix $Q_{n}\in M_{p+1} (\mathbb{R})$ 
\begin{equation}
\label{qn}
Q_n:=\begin{pmatrix}
n Id_p & -a_n\\
-a_n^T & b_n
\end{pmatrix}
\end{equation}
where  $Id_{p}$ denotes  the identity matrix in $M_{p}(\mathbb{R})$,   
\[a_n^T:=\left(\int_0^{n} \varphi_1 (t) X_t dt\dots , \int_0^{n} \varphi_p (t) X_t dt\right)\]
and
\[b_n:=\int_0^{n} X_t^2 dt.\]

Note that in the definition of the estimator $\hat{\vartheta}_n$ (\ref{tn}) stochastic integrals with respect to $X$ appear.  This integral is understood in the following sense 
\begin{equation}\label{19f-1}
\int_{0}^{t} g_{s}dX_{s}:= \int_{0} ^{t} g_{s} \left( L(s)- \alpha X_{s} \right) ds  + \int_{0} ^{t} g_{s} d Z^{H}_{s}
\end{equation}
 for every $t\geq 0$, where the second integral is a Skorohod integral with respect to the Rosenblatt process (see Section \ref{2.1}), provided that the integrals above exist.  We need to chose a Skorohod and not a pathwise integral with respect to the Rosenblatt process because, similarly to the explanation for the fBm given in e.g. \cite{HN},  the choice of the pathwise integrals (which can be easily defined for the Rosenblatt process since it has H\"older continuous paths or every order $\delta \in (0, H)$) does not lead to a consistent estimator.

First, we need to argue that the stochastic integrals that appear in (\ref{pn}) and (\ref{qn}) are well-defined. The Wiener integrals $\int_{0} ^{t} \varphi_{i} (s) dZ^{H}_{s}$  are obviously well-defined since $\varphi_{i}, i=1,.., p$ are bounded  and periodic.  In the next result we show that the Skorohod integral in (\ref{pn}) is also well-defined.

\begin{prop}
Let $(X_{t})_{t\geq 0} $ be the solution to (\ref{sde}). Then for every $t\geq 0$ the Skorohod integral $\int_{0} ^{t} X_{s} d Z^{H}_{s} $ is well-defined.
\end{prop}
{\bf Proof: }From relation  (\ref{30i-4}) in Section \ref{2.1} we need to show that 

$$
\mathbf{E} \int_{0} ^{T} \int_{0} ^{T} \Vert D_{x_{1}, x_{2}} X\Vert ^{2}_{\mathcal{H}}dx_{1} dx_{2} <\infty.
$$
 By taking the Malliavin derivative in (\ref{x}), we get for every $x_{1}, x_{2}>0$
\begin{eqnarray*}
D_{x_1x_2}X_u&=&2 d(H)1_{[0,\,u]^2}(x_1,\,x_2) I (e^{\alpha (\cdotp - u)})(x_1,\,x_2)\\
&=& 2 d(H)1_{[0,\,u]^2}(x_1,\,x_2) \int_{x_1\lor x_2}^u e^{\alpha (u'-u)}\frac{\partial K^{H'}}{\partial u'}(u',\,x_1)\frac{\partial K^{H'}}{\partial u'}(u',\,x_2)du',
\end{eqnarray*}
where $I$ is the transfer operator (\ref{I}). Hence,
\begin{eqnarray*}
 &&\|D_{x_1x_2}X\|_{\mathcal H}^2=\int_{x_1\lor x_2}^n  \int_{x_1\lor x_2}^n |u-v|^{2H-2} \\
&\times &  \int_{x_1\lor x_2}^u e^{\alpha (u'-u)}\frac{\partial K^{H'}}{\partial u'}(u',\,x_1)\frac{\partial K^{H'}}{\partial u'}(u',\,x_2)du' \int_{x_1\lor x_2}^v e^{\alpha (v'-v)}\frac{\partial K^{H'}}{\partial v'}(v',\,x_1)\frac{\partial K^{H'}}{\partial v'}(v',\,x_2)dv'dudv\\
&\leq &\int_{x_1\lor x_2}^n  \int_{x_1\lor x_2}^n |u-v|^{2H-2} \\
&\times &  \int_{x_1\lor x_2}^u \frac{\partial K^{H'}}{\partial u'}(u',\,x_1)\frac{\partial K^{H'}}{\partial u'}(u',\,x_2)du' \int_{x_1\lor x_2}^v \frac{\partial K^{H'}}{\partial v'}(v',\,x_1)\frac{\partial K^{H'}}{\partial v'}(v',\,x_2)dv'dudv\\
&=& \|D_{x_1x_2}Z^H\|_{\mathcal H}^2
\end{eqnarray*}
since $e^{\alpha (u'-u)} \leq 1$ and the other integrands are nonnegative. From Example 1 in \cite{T} we know that $\mathbf{E}  [\int_0^n\int_0^n \|D_{x_1x_2}Z\|_{\mathcal H}^2 dx_1 dx_2]<\infty$, and the result follows.
\qed

In the sequel, we will need  a more convenient expression of the estimator (\ref{tn}). Notice that the inverse of the matrix $Q_{n}$ can be expressed as (see \cite{DFW}) 
\begin{equation}
\label{30i-7}
Q_n^{-1}=\frac{1}{n}\begin{pmatrix}
Id_p + \gamma_n \Lambda_n \Lambda_n^t & -\gamma_n\Lambda_n\\
-\gamma_n \Lambda_n^t & \gamma_n
\end{pmatrix}
\end{equation}
with
\begin{equation}
\label{18f-3}
\Lambda_n  = (\Lambda_{n,\,1},\dots , \Lambda_{n,\,p})^T =\left(\frac{1}{n}\int_0^n \varphi_1 (t)X_t dt,\dots , \frac{1}{n}\int_0^n \varphi_p (t)X_t dt\right)
\end{equation}
and 
\begin{equation}
\label{18f-4}
\gamma_n = \bigg(\frac{1}{n}\int_0^n X_t^2 dt-\sum_{i=1}^p \Lambda_{n,\,i}^2\bigg)^{-1}.
\end{equation}

Another useful fact is that we can deduce a different expression for $\hat{\vartheta}_{n}$ which allows to access the error $ \hat{\vartheta}_{n}-\vartheta$ directly.

\begin{prop}\label{p2}
The estimator $\hat\vartheta_n$ (\ref{tn})  has the following representation:
\begin{equation}
\label{18f-1}\hat{\vartheta}_n = \vartheta +  Q_n^{-1}R_n
\end{equation}
with $Q_{n}$ given by (\ref{qn}) and 
\begin{equation}
\label{rn}R_n = \left(\int_0^n \varphi _1 (t) dZ^H_t,\dots , \int_0^n \varphi_p (t)dZ^H_t,\, -\int_0^n X_t dZ^H_t\right).
\end{equation}
\end{prop}
{\bf Proof: } This follows easily if the relation  $X_t=\int_0^t\left(  L(s)-\alpha X_s\right)  ds + Z^H_t$ is plugged as the integrator in each component of $P_n$ (\ref{pn}). \qed.

The relation (\ref{18f-1}) will be used in order to study the asymptotic behavior of the LSE. 

\subsection{Strong consistency}
We  study the asymptotic properties of the LSE (\ref{tn}).  In this part we prove that $\hat{\vartheta}_{n}$ is strongly consistent, i.e. it converges almost surely  to the parameter $\vartheta$ (\ref{18f-2}) as $n\to \infty$. 
In order to prove the estimator's consistency we will need several auxiliary results. First, we quote a technical lemma from \cite{KN}.
\begin{lemma} \label{l1}
Let $\gamma >0$ and $p_0\in\mathbb N$. Moreover, let $(Z_n)_{n\in\mathbb N}$ be a sequence of random variables. If for every $p\geq p_0$ there exists a constant $c_p>0$ such that for all $n\in\mathbb N$
\[(\mathbf{E}  [|Z_n|^p])^{1\slash p}\leq c_pn^{-\gamma},\]
then for all $\varepsilon >0$ there exists a random variable $\eta_\varepsilon$ such that
\[|Z_n|\leq \eta_\varepsilon n^{-\gamma + \varepsilon}\text{ a.s.}\]
for all $n\in \mathbb N$. Moreover, $\mathbf{E}  [|\eta_\varepsilon|^p]<\infty$ for all $p\geq 1$.
\end{lemma}

To show strong consistency of the estimator (\ref{tn}),   we will treat the quantities  $\frac{1}{n}R_n$ and $nQ_n^{-1}$ separately, as in  \cite{DFW} and \cite{BESV}.

\begin{prop}\label{p3}
Let $R_{n}$ be given by (\ref{rn}). Then, as $n$ tends to infinity, $\frac{1}{n}R_n\to 0$ almost surely.
\end{prop}
{\bf Proof: } Due to (\ref{30i-5}) it suffices to demonstrate  that 
\begin{equation}
\label{30i-6}\sup_n \sup_{r\in [0,\,n]}(\mathbf{E}  [|g_r|^{p^\ast}]+\mathbf{E}  [\|D^{(2)}g_r\|^{p^\ast}_{L^2([0,\,n]^2)}])<\infty
\end{equation}
for $g=\varphi_i$ ($i=1,\dots ,p$) and for $g=X$ for all $p^\ast \in\mathbb N$. Then the result will follow by taking $\gamma=1-H$ in   Lemma \ref{l1}.  Since by assumption all $\varphi_i$ are bounded, the statement for $g=\varphi_i$ ($i=1,\dots ,p$) is immediate. For $g=X$ recall that
\[X_t=\int_0^t e^{\alpha (s-t)}L(s) ds + \int_0^t e^{\alpha (s-t)}dZ^H_s.\]
Recalilng that  $L$ is bounded, we  clearly have
\[\int_0^t e^{\alpha (s-t)}L(s) ds\leq \|L\|_{\infty}\int_0^t e^{\alpha (s-t)}ds=\frac{1}{\alpha} \|L\|_{\infty} e^{-\alpha t}(e^{\alpha t}-1)\leq \frac{1}{\alpha} \|L\|_{\infty}, \]
and by the triangle inequality it is enough to prove the inequality  (\ref{30i-5}) for the random part of $X$, i.e. for $g_t=\int_0^t e^{\alpha (s-t)}dZ^H_s=Y_t $ (see (\ref{hy})). We write, for every $r>0$, 
\[\mathbf{E}  [|Y_r|^{p^\ast}]+\mathbf{E}  [\|D^{(2)}Y_r\|^{p^\ast}_{L^2([0,\,n]^2)}]=:N_{1,r}+N_{2, r}. \]

For the term $N_{1, r}$ we  note that since $Y_r$ is a multiple  Wiener-It\^o integral of order two  with respect to a Brownian motion, so  the hypercontractivity property (\ref{hypercontractivity}) is applicable, and this will give the inequality 
\[\mathbf{E}  [|Y_r|^{p^\ast}]^{1/p^{\ast}}\leq (p^\ast -1)\mathbf{E}  [|Y_r|^2]^{1/2}.\]
Therefore, since the above constant does not depend on the underlying space, it suffices to show boundedness of the $L^{2}$-norm. Due to isometry property of Wiener-Rosenblatt integrals (\ref{isoZ}) we have
\begin{eqnarray*}
\mathbf{E}  [|Y_r|^2]&=&\int_0^r\int_0^r e^{-2\alpha r} e^{\alpha u}e^{\alpha v}|u-v|^{2H-2}dudv\\
&=& \int_0^r\int_0^r e^{-\alpha u}e^{-\alpha v}|u-v|^{2H-2}dudv
\end{eqnarray*}
and clearly $\sup_{r \in [0, t] }\mathbf{E}  [|Y_r|^2]< C$ for every $t\geq 0$, with some $C>0$. 
Concerning the summand  $N_{2, r}$ we recall that
\[D_{x_1x_2}Y_r=2 d(H)1_{[0,\,r]^2}(x_1,\,x_2) I (e^{\alpha (\cdotp - r)})(x_1,\,x_2).\]
Since it is nonrandom, it is enough to prove the boundedness of $\|D^{(2)}Y_r\|^2_{L^2([0,\,n]^2)}$. We have, with $I$ given by (\ref{I}),
\begin{eqnarray*}
\|D^{(2)}Y_r\|^2_{L^2([0,\,n]^2)}&=& \int_0^n\int_0^n (2 d(H)1_{[0,\,r]^2}(x_1,\,x_2) I (e^{\alpha (\cdotp - r)})(x_1,\,x_2))^2 dx_1dx_2\\
&=&4d(H)^2\int_0^r\int_0^r ( I (e^{\alpha (\cdotp - r)})(x_1,\,x_2))^2 dx_1dx_2\\
&=&4d(H)^2 \| I (e^{\alpha (\cdotp - r)})(x_1,\,x_2)\|_{L^2([0,\,r]^2)}^2\\
&=& d(H)^2 \mathbf{E}  [I_2 ( I (e^{\alpha (\cdotp - r)})(x_1,\,x_2))^2]=d(H)^2 \mathbf{E}  [Y_r^2]
\end{eqnarray*}
due to isometry of the Wiener-It\^o integrals (\ref{isom}). As was shown above, the obtained expression is bounded by a constant independent of $r$ and of $n$. Thus, our claim (\ref{30i-6}) is proved.\qed

The next step is the almost sure convergence of the matrix $nQ_n^{-1}$. The proof is similar to the one given in \cite{DFW} for the  case of the fractional Brownian motion.

\begin{prop}\label{p4} Let $Q_{n}$ be defined by (\ref{qn}). As $n$ tends to infinity, $nQ_n^{-1}$ tends almost surely to the deterministic matrix
\begin{equation}\label{31i-1}
Q:= \begin{pmatrix} Id_{p}+ \gamma \Lambda \Lambda ^{t} & -\gamma\Lambda \\
-\gamma \Lambda ^{t}& \gamma
\end{pmatrix}  
\end{equation}
where 

\begin{equation}\label{30i-8}
\Lambda = (\Lambda_{1},..., \Lambda _{p}) \mbox{ and  } \Lambda _{i}:=\langle \varphi_{i},  \tilde{h}(t) \rangle _{ L ^{2} [0,1]}, i=1,.., p,
\end{equation}
with $\tilde{h} $ from (\ref{thy}) and
\begin{equation}\label{30i-9}
\gamma:= \left( \int_{0} ^{1}  \tilde{h}^{2}(t)dt + \alpha ^{-2H} H \Gamma (2H) -\sum_{i=1}^{p}\Lambda _{i}^2 \right) ^{-1}. 
\end{equation} 
\end{prop}
{\bf Proof: }We will use the expression (\ref{30i-7}) of the matrix $Q_{n}^{-1}$. From this formula it suffices to prove almost sure convergence of the quantities $\Lambda_{n,\,i}$ from (\ref{18f-3}) to the  constant $\Lambda _{i}$ given by (\ref{30i-8})  for every $i\in\{1,\dots , p\}$ as well as almost sure convergence of $\gamma_n^{-1}$ to the nonzero real number $\gamma^{-1}$ from (\ref{30i-9}). 
Concerning $\Lambda_{n,\,i}$ using the fact that the difference 
\[|Y_t-\tilde{Y}_t|=e^{-\alpha t}\left|\int_\mathbb R \int_\mathbb R \int_{-\infty}^0 e^{\alpha u}\frac{\partial K^{H'}}{\partial u}(u,\,x_1)\frac{\partial K^{H'}}{\partial u}(u,\,x_2)dudB(x_1)dB(x_2)\right| \]
 converges to zero almost surely as $t\to \infty$ (and the same holds true for  $|X_t-\tilde{X}_t|$), we obtain almost surely via Corollary \ref{r1}
\begin{eqnarray*}
\lim_{n\to\infty}\Lambda_{n,\,i}&=&\lim_{n\to\infty}\frac{1}{n}\int_0^n \varphi_i(t)X_tdt=\lim_{n\to\infty}\frac{1}{n}\int_0^n \varphi_i(t)\tilde{X}_tdt\\
&=&\lim_{n\to\infty}\frac{1}{n}\sum_{i=0}^{n-1}\int_{i}^{i+1}\varphi_i(t)\tilde{X}_tdt=\int_0^1 \varphi_i(t)\mathbf{E}  [\tilde X_t]dt=\int_0^1 \varphi_i(t)\tilde h(t)dt=\Lambda_{i}
\end{eqnarray*}
for every $i=1,.., p$. Concerning  $\gamma_n^{-1}$ we have from (\ref{18f-4})
\begin{eqnarray*}
\frac{1}{n}\int_0^n X_t^2 dt &=&\frac{1}{n}\int_0^n h(t)^2dt+\frac{2}{n}\int_0^n h(s)Y_s ds+\frac{1}{n}\int_0^n Y_s^2ds.
\end{eqnarray*}
Since $|h(t)-\tilde{h}(t)|=e^{-\alpha t}|\int_{-\infty}^0 e^{\alpha s}L(s)ds|$, we conclude that the first integral converges to $\int_0^1 \tilde{h}^2(t)dt$. For the second integral note that due to boundedness of $\frac{1}{n}\int_0^n Y_s ds$ (shown in \cite{ntra}) and of $|\frac{1}{n}\int_0^n h(t)dt|$ we obtain almost surely
\[\lim_{n\to\infty}\frac{2}{n}\int_0^n h(s)Y_s ds=\lim_{n\to\infty}\frac{2}{n}\int_0^n \tilde h(s)\tilde Y_s ds=0\]
by applying  Proposition \ref{pp1}. The almost sure limit of the third integral equals $\alpha^{-2H}H\Gamma (2H)$, as demonstrated in \cite{ntra}.
So  almost surely
 $$\gamma_n^{-1}=\frac{1}{n}\int_0^n X_t^2 dt-\sum_{i=1}^p \Lambda_{n,\,i}^2\to _{n\to \infty}  \|\tilde{h}\|_{L^2([0,\,1])}-\sum_{i=1}^p\langle \tilde h,\,\varphi_i \rangle ^2_{L^2([0,\,1])}+\alpha^{-2H}H\Gamma (2H)$$ and by Bessel's inequality we can see as in \cite{DFW} that the above limit  is indeed a positive real number.\qed

As a consequence of Propositions \ref{p2}, \ref{p3} and \ref{p4} we obtain the strong consistency of the least squares estimator.

\begin{theorem}
As $n\to \infty$, the LSE (\ref{tn}) converges almost surely to the parameter $\vartheta =\left( \mu _{1},..., \mu _{p}, \alpha\right).$ 
\end{theorem}

\section{ Limit distribution of the least squares estimator}\label{5}
We will analyze the asymptotic behavior in distribution of the LSE. We use the decomposition  of $\hat{\vartheta}_{n}$ given in Proposition \ref{p2}. It follows from this result, since the random matrix  $n Q_{n} ^{-1}$ given by (\ref{qn}) converges almost surely to the deterministic matrix $Q$ from Proposition \ref{p4},  it is enough to consider the asymptotics of the vector $R_n$ in (\ref{rn}).

We start with a result concerning the first $p$ components of the vector (\ref{rn}). In the sequel, by a Rosenblatt random variable we mean a random variable with the same law as $ Z^{H}_{1}$ from (\ref{rose}).

\begin{prop}\label{p6}
For every $n\geq 1$, consider $U_{n}: =n^{-H}\int_0^n f(s)dZ^H_s$ for a bounded $1$-periodic function $f$. As $n$ tends to infinity, this sequence converges in distribution to $U=\left(\int_0^1 f(t)dt\right)  V$, where $V$ is a Rosenblatt random variable.
\end{prop}
{\bf Proof: }
It follows by the scaling property of the Rosenblatt process that $U_n\stackrel{d}{\equiv}\int_0^n f(ns)dZ^H_s$, where $\stackrel{d}{\equiv}$ stands for the equivalence of finite dimensional distributions. We will show that this sequence converges  in $L^2 (\Omega)$ to the random variable $\left(\int_0^1 f(t)dt\right) Z^H_1$. We can write 

\begin{eqnarray*}
&& \mathbf{E} \left[ \int_{0} ^{1} f(ns) dZ^{H} _{s}-\left(\int_{0} ^{1} f(s)ds \right) Z^{H}_{1}\right] ^{2} = \mathbf{E}\left[ \int_{0} ^{1}\left( f(ns)- \int_{0} ^{1} f(r)dr \right) dZ^{H}_{s} \right] ^{2} \\
&=&H(2H-1)  \int_{0} ^{1} \int_{0} ^{1} dudvf(nu)f(nv) \vert u-v\vert ^{2H-2} +\left(\int_{0} ^{1} f(s)ds\right) ^{2} \\ &&-2H(2H-1) \int_{0} ^{1} \int_{0} ^{1}dudv  f(nu) \vert u-v\vert ^{2H-2} \int_{0} ^{1} f(s)ds.
\end{eqnarray*} 

\b First, 

\begin{eqnarray*}
&&H(2H-1)  \int_{0} ^{1} \int_{0} ^{1} dudvf(nu)f(nv) \vert u-v\vert ^{2H-2}  =H(2H-1) n ^{-2H} \int_{0} ^{n} \int_{0} ^{n} dudv  f(nu) f(nv) \vert u-v \vert ^{2H-2} \\
&=& n ^{-2H} H(2H-1) \sum_{i,j=0} ^{n-1}  \int_{0} ^{1} \int_{0} ^{1} dudvf(u)f(v) \vert u-v +i-j\vert ^{2H-2}\\
&\sim&  n ^{-2H} H(2H-1) \sum_{i,j=0,\,i\neq j} ^{n-1}  \int_{0} ^{1} \int_{0} ^{1} dudvf(u)f(v) \vert i-j\vert ^{2H-2}\left\vert 1 + \frac{u-v}{i-j}\right\vert ^{2H-2}\\
&\sim & n ^{-2H} H(2H-1) \sum_{i,j=0; i\not=j} ^{n-1} \vert i-j\vert ^{2H-2}  \left(\int_{0} ^{1} f(s)ds\right) ^{2} \to _{n\to \infty}  \left(\int_{0} ^{1} f(s)ds\right) ^{2} ,
\end{eqnarray*}
the symbol $\sim$ signifying asymptotic equivalence, i.e., both sides having the same limit as $n$ tends to infinity. The equivalence is obtained by considering the binomial expansion of $\left\vert 1 + \frac{u-v}{i-j}\right\vert ^{2H-2}$. On the other hand,

\begin{eqnarray*}
&&H(2H-1) \int_{0} ^{1} \int_{0} ^{1}dudv  f(nu) \vert u-v\vert ^{2H-2}  \\
&=& H(2H-1) \left( \int_{0}^{1} f(nu)du \int_{0} ^{u} (u-v) ^{2H-2} dv + \int_{0}^{1} \int_v^1 f(nu) du (v-u)^{2H-2} dv \right) 
\\
&=& H \int_{0}^{1} f(nu) du u ^{2H-1}+ H\int_{0} ^{1} f(nu)du (1-u) ^{2H-1}.
\end{eqnarray*}

Now, again by the binomial expansion,
\begin{eqnarray*}
&&H \int_{0}^{1} f(nu) du u ^{2H-1} =H n ^{-2H} \int_{0} ^{n} f(u) u ^{2H-1} du \\
&=& H n ^{-2H} \sum_{i=0} ^{n-1} \int_{0} ^{1} f(u) (u+i) ^{2H-1} du \sim H n ^{-2H} \int_{0} ^{1} f(u) \sum_{ i=0} ^{n-1} (u+i) ^{2H-1} du \\
&\sim  &  H n ^{-2H} \int_{0} ^{1} f(u) \sum_{ i=1} ^{n-1} i^{2H-1}\left(1+\frac{u}{i}\right) ^{2H-1} du \sim H n ^{-2H}  \int_{0} ^{1}f(u) du\frac{n^{2H}}{2H} \\
&=& \frac{1}{2} \int_{0}^{1} f(u)du.
\end{eqnarray*}

Moreover,
\begin{eqnarray*}
 &&H\int_{0} ^{1} f(nu)du (1-u) ^{2H-1}=H\int_0^1 f(n(1-u)) u^{2H-1} du\\
&=& H\int_0^1 f(-nu) u^{2H-1} du \to _{n\to \infty}\frac{1}{2} \int_{0}^{1} f(-u)du=\frac{1}{2} \int_{0}^{1} f(u)du
\end{eqnarray*}
with the same argument as above. This gives the desired  $L^{2}(\Omega)$-convergence.
\qed

Now let us consider the last component of the vector $R_n$ in (\ref{rn}). First we show that the stochastic integral  part  does not contribute to the limit.

\begin{prop}\label{pp7}Let $(Y_{t})_{t\geq 0}$ be given by (\ref{hy}). Then, as $n$ tends to infinity, 
$$\mathbf{E}  \left(n^{-H}\int_0^n Y_t dZ^H_t\right)^2\to 0.$$
\end{prop}
{\bf Proof: } Let us estimate the $L^{2}$-norm of the random variable $n^{-H}\int_0^n Y_t dZ^H_t$ with $Y$ from (\ref{hy}).  In \cite{T} the following bound is given:
\begin{eqnarray*}
\mathbf{E}   (\int_0^n Y_t dZ^H_t)^2& \leq& C \bigg( \mathbf{E} \left[\int_0^n \int_0^n Y_u Y_v |u-v|^{2H-2}du dv\right]\\
&&\qquad + \mathbf{E} \left[\int_0^n\int_0^n \int_0^n \int_0^n D_{x_1,\, x_2}Y_u   D_{x_1,\, x_2}Y_v |u-v|^{2H-2}du dv dx_1 dx_2\right]\bigg).
\end{eqnarray*}
Since $Y_{u}$ is a double integral, it is easy to note that  the two summands  above only differ by a constant, so it is enough to consider one of them. We obtain using the isometry for the Rosenblatt process
\begin{eqnarray*}
&&\mathbf{E} \left[\int_0^n \int_0^n Y_u Y_v |u-v|^{2H-2}du dv\right]\\
&=& \int_0^n\int_0^n \int_0^u\int_0^v e^{\alpha (s-u)}e^{\alpha (r-v)}|r-s|^{2H-2}dr ds |u-v|^{2H-2}du dv\\
&\leq & \int_0^n\int_0^n \int_0^n\int_0^n e^{-\alpha |s-u|}e^{-\alpha |r-v|}|r-s|^{2H-2}|u-v|^{2H-2}dr ds du dv,
\end{eqnarray*}
and it was demonstrated in \cite{HN} and \cite{HNZ} that this bound multiplied by $\frac{1}{n}$, $\frac{1}{n\log (n)}$ or $n^{2-4H}$ in cases $H\in \left(\frac{1}{2},\,\frac{3}{4}\right)$, $H=\frac{3}{4}$ and $H>\frac{3}{4}$ respectively converges to a constant. Thus, the statement follows.\qed

The next proposition concludes the asymptotic analysis.

\begin{prop}Let  $(Y_{t})_{t\geq 0}$ be given by (\ref{x}).
The sequence $n^{-H}\int_0^n X_t dZ^H_t$ converges in distribution to $U=\left(\int_0^1 \tilde{h}(t)dt\right)  V$, where $V$ is a Rosenblatt random variable.
\end{prop}
{\bf Proof: } Recall that for every $t\geq 0$, $X_{t}= Y_{t}+ h(t)$, see (\ref{hy}), so we need to analyze the limit of $n^{-H}\int_0^n  h(t)  dZ^H_t$. Since $\tilde{h}$ from (\ref{thy}) is a periodic function, it suffices to demonstrate that $n^{-H}\int_0^n (h(t)-\tilde{h}(t)) dZ^H_t$ converges to zero in $L^2 (\Omega)$ and then to apply Proposition \ref{p6}. Since $|h(t)-\tilde{h}(t)|$ is bounded by $e^{-\alpha t}$ times a constant, we get by the isometry property (\ref{isom})
\[\mathbf{E}  \left[\left(\int_0^n (h(t)-\tilde{h}(t)) dZ^H_t\right)^2\right]\leq c\int_0^n \int_0^n e^{-\alpha u}e^{-\alpha v}|u-v|^{2H-2}du dv,\]
which is bounded uniformly in $n$, and the desired  convergence follows. \qed

By putting together the above results, we state and prove the main result of this section.

\begin{theorem}
Let $\hat{\vartheta }_{n}$ be given by (\ref{tn}).  Then the sequence $n^{1-H}\left( \hat{\vartheta }_{n} -\vartheta \right)$
converges in distribution, as $n\to \infty$, to $QR$ where the matrix $Q$ is given by (\ref{31i-1}) and $R$ is the following random vector
$$R= \left( \int_{0} ^{1} \varphi_{1}(s)ds,...,  \int_{0} ^{1} \varphi_{p}(s)ds, -\int_{0} ^{1}\tilde{h}_{s}ds  \right)^T V$$
where $V$ is a Rosenblatt ranfom variable and $\tilde{h}$ is defined by (\ref{thy}). 
\end{theorem}
{\bf Proof: } The almost sure convergence of $nQ_{n} ^{-1}$ to the matrix $Q$ follows from  Proposition \ref{p4} an we need to prove the asymptotic behavior in distribution of the vector $\frac{1}{n}R_{n}$ (\ref{rn}).  For any $a_1,\dots , a_{p+1}\in \mathbb R$ and for $1$-periodic functions $f_1,\dots ,f_{p+1}$ we have
\[\sum_{i=1}^{p+1} a_i n^{-H}\int_0^n  f_i(t) dZ^H_t=  n^{-H}\int_0^n  \sum_{i=1}^{p+1} a_i f_i(t) dZ^H_t,\]
and by Proposition \ref{p6} this converges in distribution as $n\to \infty$ to $U=(\int_0^1  \sum_{i=1}^{p+1} a_i f_i(t) dt)V$ (where $V$ is a Rosenblatt random variable), because $ \sum_{i=1}^{p+1} a_i f_i$ is again a $1$-periodic function. By applying the results to $f_{i}=\varphi_{i}, i=1,.., p$ and $f_{p+1}= -\tilde{h}$ and by using the $L^ {2}$ convergence from Proposition \ref{pp7}, we obtain the conclusion. \qed

Note that for functions $\varphi_i$, $i=1,\dots ,p$, whose integrals are equal to zero one might obtain an improvement in the speed of convergence. This case is, however, not treated here.

\section{Different estimators}
The estimator $\vartheta _{n}$ (\ref{tn}), although consistent and with explicit limit distribution, involves a Skorohod integral. It is well-known that it is difficult to simulate such a stochastic object. Therefore, we will define some alternative estimators that can be expressed only in terms of Wiener and Lebesque integrals and consequently they can be simulated. One of these new estimators represents an extended version of the estimators proposed in \cite{HN} or \cite{ntra} as it reduces to them when the periodic drift $L$ reduces to a constant.  

Recall that the the functions $\varphi_{i}$ from (\ref{L}) are assumed to be orthogonal in $L^ {2} ([0,1])$. We will consider the following assumptions (the function $\tilde{h} $ is defined in (\ref{thy})):
\begin{itemize}
\item[(A1)] $\tilde{h} $ does not belong to $\operatorname{span}(\varphi_1,\dots ,\varphi_p)$. In this case there exists a bounded function $\varphi_{p+1}$ orthogonal to all $\varphi_i$ ($i\in \{1,\,p\}$), but not orthogonal to $\tilde{h}$.
\item[(A1*)]  $\tilde{h}\in \operatorname{span}(\varphi_1,\dots ,\varphi_p)$. Then there  is no $L^2$ function satisfying the above orthogonality conditions.
\end{itemize}
We will show below in Remark \ref{r2} that in the case when $\varphi_{i}, i=1,..,p$ are elements of the trigonometric basis of $L^ {2} ([0, 1])$, then it is easy to check which of these assumptions is satisfied and to determine the function $\varphi_{p+1}$ without the knowledge of $\tilde{h}$ in case of (A1).

\begin{prop}\label{consA1}
Assume that (A1) is satisfied. Define for every $n\geq 1$ 
\[\bar{\alpha}_{n} :=-\frac{\int_0^n \varphi_{p+1}(t)dX_t}{\int_0^t \varphi_{p+1}(t)X_t dt}\]
and for $i=1,.., p$
\begin{equation*}
\bar{\mu}_{i, n}:=\frac{1}{n} \left(\int_0^n \varphi_{i}(t)dX_t + \bar{\alpha}\int_0^n \varphi_{i}(t)X_t dt\right).
\end{equation*}
Then $ \left(\bar{\alpha}_{n}, \bar{\mu} _{1, n}, ..., \bar{\mu}_{p, n} \right)$ is a 
 consistent estimator of the parameter  $(\alpha, \mu_{1},.., \mu _{p}) $ of  the model (\ref{sde}). 
\end{prop}
{\bf Proof: } From (\ref{19f-1}) and (A1) we have 
\[\frac{1}{n}\int_0^n \varphi_{p+1}(t)dX_t=-\alpha \frac{1}{n}\int_0^n \varphi_{p+1}(t)X_t dt+\frac{1}{n}\int_0^n \varphi_{p+1}(t)dZ^H_t\]
so we can write 
\begin{equation}\label{19f-2}
\bar{\alpha}_{n} -\alpha= \frac{ n ^ {-1} \int_{0} ^ {n} \varphi_{p+1} (t) dZ ^ {H} _{t} }{ n ^ {-1} \int_{0} ^ {n} \varphi_{p+1} (t) X_{t} dt}.
\end{equation}

As demonstrated in Proposition \ref{p3}, the numerator of (\ref{19f-2})  converges to zero almost surely as $n\to \infty$. Moreover, we can conclude using Proposition \ref{pp1} that
\[\Lambda_{n,\,p+1}:=\frac{1}{n}\int_0^n \varphi_{p+1}(t)X_t dt\to _{n \to \infty} \langle \tilde{h},\, \varphi_{p+1}\rangle_{L^2([0,\,1])}\]
almost surely. Since this is nonzero by the assumption (A1), strong consistency of $\bar{\alpha}_{n}$ follows. Consistency of $\bar{\mu}_i$ follows by observing that
\[\frac{1}{n}\int_0^n \varphi_{i}(t)dX_t=\mu_i- \alpha \frac{1}{n}\int_0^n \varphi_{i}(t)X_t dt+\frac{1}{n}\int_0^n \varphi_{i}(t)dZ^H_t,\]
and this implies, for every $i=1,.., p$
\begin{equation}\label{19f-3}
\bar{\mu }_{i, n}-\mu _{i}=\frac{1}{n}(\bar{\alpha}_{n} -\alpha) \int_{0} ^ {n}  \varphi_{i}(t) X_{t}dt + \frac{1}{n} \int_0^n \varphi_{i}(t)dZ^H_t
\end{equation}
and  the last summand again converges to zero almost surely as $n\to \infty$ while $\frac{1}{n} \int_{0} ^ {n}  \varphi_{i}(t) X_{t}dt$ tends to a constant.\qed

The asymptotic behavior in distribution of the above estimators can be easily obtained from the proofs in Section \ref{5}. 
\begin{prop}
As $n$ tends to infinity the vector $n^{1-H}(\bar{\alpha}_{n}  -\alpha ,\,\bar{\mu}_{1,n}  -\mu_1,\dots , \bar{\mu}_{p, n}-\mu_p )^T$ converges in distribution to the vector
\[ \begin{pmatrix}\int_0^1 \varphi_{p+1}(t)dt \frac{1}{\langle \varphi_ {p+1},\, \tilde h \rangle _{L^2([0,\,1])}}\\ \int_0^1 \varphi_{p+1}(t)dt \frac{\langle \varphi_ {1},\, \tilde h \rangle _{L^2([0,\,1])}}{\langle \varphi_ {p+1},\, \tilde h \rangle _{L^2([0,\,1])}}+\int_0^1 \varphi_1 (t) dt\\ \vdots \\ \int_0^1 \varphi_{p+1}(t)dt \frac{\langle \varphi_ {p},\, \tilde h \rangle _{L^2([0,\,1])}}{\langle \varphi_ {p+1},\, \tilde h \rangle _{L^2([0,\,1])}}+\int_0^1 \varphi_p (t) dt   \end{pmatrix}V, \]
where $V$ is a Rosenblatt random variable.
\end{prop}
{\bf Proof: }
This follows by construction from relations (\ref{19f-2}), (\ref{19f-3}),  Proposition \ref{pp1} and the non-central limit theorem in Proposition \ref{p6}.\qed

When the assumption (A1*) is satisfied, we can also define consistent estimators for the parameters of the model (\ref{sde}) which involve only Wiener and deterministic integrals.
\begin{prop}\label{pp11}
Assume that (A1*) is satisfied. Consider the following estimators 
\[\bar{\alpha}_{n} ^{(1)}:=\left(\frac{1}{H\Gamma (2H)}\gamma_n^{-1}\right)^{-\frac{1}{2H}}\]
and for $i=1,.., p$,
\[\bar{\mu}_{n,i}^{(1)}:=\frac{1}{n} \left(\int_0^n \varphi_{i}(t)dX_t+\bar{\alpha}^{(1)} \int_0^t \varphi_{i}(t)X_t dt\right)\]
Then $\left( \bar{\alpha}_{n} ^{(1)}, \bar{\mu}_{1, n} ^{(1)},...,  \bar{\mu}_{p,n} ^{(1)}\right)$  is a  consistent estimator of the parameter (\ref{18f-2}).
\end{prop}
{\bf Proof: }
It was shown in Proposition \ref{p4} that with $\gamma_{n}$ defined in (\ref{18f-4})
\[\gamma_n^{-1}\to _{n\to \infty}\|\tilde{h}\|_{L^2([0,\,1])}-\sum_{i=1}^p\langle \tilde h,\,\varphi_i \rangle ^2_{L^2([0,\,1])}+\alpha^{-2H}H\Gamma (2H)\]
almost surely. Because (A1*) is satisfied, we obtain the equality $ \|\tilde{h}\|_{L^2([0,\,1])}=\sum_{i=1}^p\langle \tilde h,\,\varphi_i \rangle ^2_{L^2([0,\,1])}$, and thus consistency follows by the continuous mapping theorem. Consistency of the estimators of the $\mu_i$ is a direct consequence and can be shown similarly to the strong consistency in Proposition \ref{consA1}.\qed

Concerning the limit in law of $\left( \bar{\alpha}_{n} ^{(1)}, \bar{\mu}_{1, n} ^{(1)},...,  \bar{\mu}_{p,n} ^{(1)}\right)$, we have the following result.
\begin{prop}
As $n$ tends to infinity the vector $n^{1-H}(\bar{\alpha}_{n}^{(1)}-\alpha ,\,\bar{\mu}^{(1)}_{1, n}-\mu_1,\dots , \bar{\mu}^{(1)}_{p, n}-\mu_p )^T$ converges in distribution to the vector
\[C_{\alpha}G_{\infty } \begin{pmatrix}1\\ \langle \tilde h ,\, \varphi_1 \rangle_{L^2([0,\,1])} \\ \vdots \\ \langle \tilde h ,\, \varphi_p \rangle_{L^2([0,\,1])}   \end{pmatrix}+ Z^H_1 \begin{pmatrix} 0 \\ \int_0^1 \varphi_1 (t)dt \\ \vdots \\ \int_0^1 \varphi_p (t)dt \end{pmatrix}, \]
where $C_\alpha = \frac{\alpha^H}{2H^2 \Gamma (2H)}$ and $G_\infty = B_H \times R$ with $R$ being $\sigma (Z^H)$-measurable and having a Rosenblatt distribution and $B_H$ being defined as follows:
\[B_H=\frac{(2H-1)\Gamma (H+1)}{\sqrt{\frac{H}{2}(2H-1)}}.\]
\end{prop}
{\bf Proof: } Using a Taylor expansion we obtain for large $n$
\begin{eqnarray*}
\bar{\alpha}_{n}^{(1)}-\alpha&=&\alpha \left(\left(1+\frac{\alpha^{2H}(\gamma_n^{-1}-\alpha^{-2H}H\Gamma (2H))}{H \Gamma (2H)}\right)^{-\frac{1}{2H}}-1\right)\\
&=&\frac{\alpha^{2H+1}}{2H^2 \Gamma (2H)}(\gamma_n^{-1}-\alpha^{-2H}H\Gamma (2H))+o(1).
\end{eqnarray*}
Therefore, it suffices to calculate the asymptotics of the quantity
\[\gamma_n^{-1}-\alpha^{-2H}H\Gamma (2H)=\frac{1}{n}\int_0^n X_t^2 dt -\sum_{i=1}^p \left(\frac{1}{n}\int_0^n X_t \varphi_i(t)dt\right)^2 -\alpha^{-2H}H\Gamma (2H).\]
As in the previous computations, the above expression has the same limit in distribution, as $n \to \infty$, as
\begin{eqnarray}
&&\left(\frac{1}{n}\int_0^n \tilde X_t^2 dt -\sum_{i=1}^p \left(\frac{1}{n}\int_0^n \tilde X_t \varphi_i(t)dt\right)^2 -\alpha^{-2H}H\Gamma (2H)\right)\nonumber \\
&=&\frac{1}{n}\int_0^n \tilde{Y}_t^2dt-\frac{2}{n}\int_0^n \tilde{Y}_t\tilde h(t)dt+\frac{1}{n}\int_0^n \tilde{h}(t)^2dt-\alpha^{-2H}H\Gamma (2H)\nonumber \\
&&\qquad -\sum_{i=1}^p \left(\frac{1}{n}\int_0^n \tilde Y_t \varphi_i(t)dt\right)^2+2\sum_{i=1}^p \left(\frac{1}{n}\int_0^n \tilde Y_t \varphi_i(t)dt\right)\left(\frac{1}{n}\int_0^n \tilde h(t) \varphi_i(t)dt\right)-\sum_{i=1}^p \langle\tilde h,\, \varphi_i\rangle_{L^2([0,\,1])}^2 \nonumber \\
&=&\frac{1}{n}\int_0^n \tilde{Y}_t^2dt-\frac{2}{n}\int_0^n \tilde{Y}_t\tilde h(t)dt-\alpha^{-2H}H\Gamma (2H)\nonumber \\
&&\qquad -\sum_{i=1}^p \left(\frac{1}{n}\int_0^n \tilde Y_t \varphi_i(t)dt\right)^2+2\sum_{i=1}^p \left(\frac{1}{n}\int_0^n \tilde Y_t \varphi_i(t)dt\right)\left(\frac{1}{n}\int_0^n \tilde h(t) \varphi_i(t)dt\right).\label{19f-4}
\end{eqnarray}
Note that $\frac{1}{n}\int_0^n \tilde{h}(t)^2dt$ and $\sum_{i=1}^p \langle\tilde h,\, \varphi_i\rangle_{L^2([0,\,1])}^2$ cancel each other out by Parseval's identity due to (A1*). If we consider the space $[0,\,n]$ with the scalar product
\[\langle f,\, g \rangle_n:=\frac{1}{n}\int_0^n f(x)g(x)dx ,\]
the orthonormality assumption of $\varphi_{i}$,  as well as (A1*),  will still hold for the periodic extensions on $[0, n]$ of $\varphi_i$ and $\tilde h$ under the scalar product $\langle \cdot, \cdot\rangle_{n}$, and  by the assumption (A1*) we obtain
\[2\sum_{i=1}^p \left(\frac{1}{n}\int_0^n \tilde Y_t \varphi_i(t)dt\right)\left(\frac{1}{n}\int_0^n \tilde h(t) \varphi_i(t)dt\right)=2\langle \tilde h,\,\tilde Y\rangle_n=\frac{2}{n} \int_0^n \tilde{Y}_t \tilde h(t)dt.\]
Therefore, (\ref{19f-4}) reduces to the term
\begin{eqnarray*}
&&\frac{1}{n}\int_0^n \tilde{Y}_t^2dt-\alpha^{-2H}H\Gamma (2H)-\sum_{i=1}^p \left(\frac{1}{n}\int_0^n \tilde Y_t \varphi_i(t)dt\right)^2\\
&=&\frac{1}{n}\int_0^n \left( \tilde{Y}_t^2-\mathbf{E} [\tilde Y_t^2]\right) dt+\frac{1}{n}\int_0^n\left( \mathbf{E} [\tilde{Y}_t^2]-\alpha^{-2H}H\Gamma (2H) \right)dt-\sum_{i=1}^p \left(\frac{1}{n}\int_0^n \tilde Y_t \varphi_i(t)dt\right)^2.
\end{eqnarray*}
It follows from Proposition \ref{pp1} that $n^{1-H}\sum_{i=1}^p \left(\frac{1}{n}\int_0^n \tilde Y_t \varphi_i(t)dt\right)^2$ converges to zero in $L^2(\Omega) $ an $n\to \infty$. As to the first two summands, by replacing once again $\tilde{Y}$ by $Y$, the quantity (\ref{19f-4}) will become asymptotically equivalent to
\[\frac{1}{n}\int_0^n \left(  {Y}_t^2-\mathbf{E} [Y_t^2]\right) dt+\frac{1}{n}\int_0^n \left(\mathbf{E} [{Y}_t^2]-\alpha^{-2H}H\Gamma (2H)\right) dt.\]
It has been shown in \cite{ntra} that $n^{1-H}\frac{1}{n}\int_0^n \mathbf{E}\left(  [{Y}_t^2]-\alpha^{-2H}H\Gamma (2H) \right) dt$ goes to zero in $L^2 (\Omega)$ when $n\to \infty$. Another result from \cite{ntra} by rescaling of $Z^H$ by the factor $n^{-H}$ is that
\[n^{1-H}\frac{1}{n}\int_0^n \left( {Y}_t^2-\mathbf{E} [Y_t^2]\right) dt \stackrel{d}{\equiv}\alpha ^{-H-1}G_{\alpha n}\]
where $G_T$ are explicitly defined random variables converging in $L^2$ as $T\to \infty$ to a limit denoted by $G_{\infty}$, whose distribution and properties are as claimed in the statement of the  proposition. Thus, as $n\to \infty$
\[ n^{1-H}\frac{\alpha^{2H+1}}{2H^2 \Gamma (2H)} (\bar{\alpha}_{n}^{(1)}-\alpha)\stackrel{d}{\to}\alpha ^{-H-1} G_{\infty}.\]

By the definition of $\bar{\mu}_i^{(1)}$, we can write for every $i=1,.., p$
\begin{equation*}
\bar{\mu}_{i,n} ^ {(1)} -\mu_{i}= (\bar{\alpha}_{n} ^ {(1)} -\alpha) \frac{1}{n} \int_{0} ^ {n} \varphi_{i}(t) X_{t}dt + \frac{1}{n} \int_{0} ^ {n} \varphi_{i} (t) dZ^ {H} (t).
\end{equation*}

 Since  the sequence $\frac{1}{n} \int_{0} ^ {n} \varphi_{i}(t) X_{t}dt$ converges almost surely as $n\to \infty$ to $ \langle \tilde h,\, \varphi_i \rangle_{L^2([0,\,1])}$,  it now suffices to investigate joint convergence of 
\begin{equation*}
\left( \frac{1}{n}\int_0^n 
\left( {Y}_t^2-\mathbf{E} [Y_t^2]\right) dt, \frac{1}{n}\int_0^n f(s)dZ^H_s\right)
\end{equation*} for a periodic function $f$. First we rescale the Rosenblatt process involved in both elements by $n^{-H}$ and obtain
\begin{equation}
\label{19f-6}\left(n^{1-H}\frac{1}{n}\int_0^n \left(  {Y}_t^2-\mathbf{E} [Y_t^2]\right) dt ,\, n^{-H}\int_0^n f(s) dZ^H_s \right)\stackrel{d}{\equiv}\left(\alpha ^{-H-1}G_{\alpha n} ,\,  \int_{0} ^{1} f(ns) dZ^{H}_{s}  \right).
\end{equation}

We know from Proposition \ref{p6} that $\int_{0} ^{1} f(ns) dZ^{H}_{s}$ converges in $L^2$ to $(\int_0^1 f(s)ds)Z^H_1$, and the first component also converges in $ L ^ {2}$, as mentioned above. Consequently, we get   the joint convergence in distribution  of the vector (\ref{19f-6}) to $(\alpha ^{-H-1}G_{\infty},\,(\int_0^1 f(s)ds)Z^H_1 )$. This fact  combined with Slutsky's lemma for vectors yields the desired result.
\qed

The random vector $ (G_{\infty}, Z^{H}_{1})  $ whose components appear in the statement of the above result can be understood as a two dimensional Rosenblatt vector. Its marginals are Rosenballt distributed and it is well-defined as a limit in $ L ^{2}(\Omega)$ of the sequence (\ref{19f-6}).

Let us end by a discussion concerning the hypotheses (A1) and (A1*) in the case of the trigonometric basis of $L^ {2}([0,1]).$

\begin{remark}\label{r2}
\begin{itemize}
\item Consider  the orthonormal  basis of $L^ {2}([0,1])$ formed by  \\
$\{1,\,\sqrt 2 \sin (2\pi n \cdot),\,\sqrt 2 \cos (2\pi n \cdot),\, n\in\mathbb N\}$. Recall that $\tilde{h}(t)= \sum_{i=1} ^ {p}\mu _{i} \int_{0} ^ {t} e ^ {-\alpha (t-s)} \varphi_{i}(s)ds.$ By direct calculation, we obtain

\begin{eqnarray*}
\int_{-\infty}^t e^{\alpha (s-t)}\sin (2\pi n s)ds = \frac{\alpha}{(2\pi n)^2 +\alpha^2}\sin (2\pi n t)-\frac{2\pi n}{(2\pi n)^2+\alpha^2}\cos (2\pi n t),\\
\int_{-\infty}^t e^{\alpha (s-t)}\cos (2\pi n s)ds = \frac{\alpha}{(2\pi n)^2 +\alpha^2}\cos (2\pi n t)+\frac{2\pi n}{(2\pi n)^2+\alpha^2}\sin (2\pi n t).
\end{eqnarray*}
This implies a simple rule: If $\{\varphi_1,\dots ,\varphi_p\}$ are elements of the trigonometric basis and if this set is ''symmetric'' (i.e., $\sin (2\pi n \cdot)\in \{\varphi_1,\dots ,\varphi_p\} \Leftrightarrow \cos (2\pi n \cdot)\in \{\varphi_1,\dots ,\varphi_p\} $), then the assumption (A1*) is satisfied; otherwise, (A1) is verified and $\varphi_{p+1}$ can be chosen from the missing counterparts.

\item The pathwise estimators of $\alpha$ considered in \cite{HN} and \cite{ntra} are special cases of the  estimator defined in Proposition \ref{pp11}. Indeed, for a constant mean function the assumption (A1*) is satisfied.
\end{itemize}
\end{remark}

\section{Appendix: The basics of the Malliavin calculus}

Here we present the tools from Malliavin calculus needed throughout the paper. See \cite{N} or \cite{NPbook} for more details.

\subsection{Multiple Wiener-It\^o Integrals}\label{sec21}
Let $(B_{t})_{t\in [0, T]}$ be a Brownian motion defined on a probability space $\left( \Omega ,\mathcal{F},\mathbb{P}\right) $. For a deterministic function $h\in
L^{2}\left( [0,T]\right)$, the Wiener integral $
\int_{0}^{T} h\left( s\right) dB\left( s\right) $ is also denoted by
$B(h) $. The inner product $\int_{0}^{T} f\left(
s\right) g\left( s\right) ds$ will be denoted by $\left\langle
f,g\right\rangle_{L^{2}\left([0,T]\right)}$.
\\
For every $q\geq 1$, ${\mathcal{H}
}^{B}_{q}$ denotes the $q$th Wiener chaos of $B$, defined as the closed linear
    subspace of $L^{2}(\Omega )$ generated by the random variables $%
    \{H_{q}(B(h)),h\in L ^{2}([0, T]),\Vert h\Vert _{L ^{2} ([0,T]^{q})}=1\}$ where $%
    H_{q}$ is the $q$th Hermite polynomial. \\
    The mapping ${I_{q}(h^{\otimes q}):}=q!H_{q}(B(h))$ can be extended to a linear isometry between $L^{2}_{s}(\mathbb{R}^{q})$ the space of symmetric square integrable functions of $[0,T]^{q}$ (equipped
    with the modified norm ${
        \sqrt{q!}}\Vert .\Vert _{L^{2}([0,T]^{q})}$) and $\mathcal{H}^{B}_{q}$%
    . When $f \in L^{2}([0,T]^{q})$, the random variable ${%
        I_{q}(f)}$ can be interpreted as a multiple Wiener-It\^o integral of
    ${f}$ of order $q$ w.r.t. $B$ and in this case, we write :
    \begin{equation}\label{IntegralMultiple}
    I_{q}(f) = \int_{[0,T]^{q}} f(y_{1},\ldots,y_{q})
    dB _{y_{1}}\ldots dB_{y_{q}}.
    \end{equation}
From the many properties of multiple Wiener-It\^o integrals  we recall now two that we will need in our study. The first one is the \textbf{isometry property}, which states that for every $f \in L^{2}_{s}([0,T]^{q})$, $g\in
L^{2}_{s}([0,T]^{p})$, with $p,q \geq 1$, the following holds:
    \begin{equation}\label{isom}
    \mathbf{E}\left[ I_{q}(f)I_{p}(g)\right] =
    \begin{cases}
     p! \times \big\langle f, g \big\rangle_{L^2([0,T]^{p})} \qquad &\text{if } p=q,\\
    0 \qquad\qquad\qquad\quad\quad &\text{if } p \ne q.
    \end{cases}
    \end{equation}%
The second one is the \textbf{hypercontractivity property} which states that for $f \in L^{2}_{s}([0,T]^{q})$, $q\geq1$,
 the multiple Wiener-It\^o integral $I_{q}(f)$ satisfies a hypercontractivity property (equivalence
    in $\mathcal{H}^{B}_{q}$ of all $L^{p}(\Omega)$ norms for all $p\geq 2$), which
    implies that for any $F\in \oplus _{l=1}^{q}\mathcal{H}^{B}_{l}$ (i.e. in a
    fixed sum of Wiener chaoses), we have
    \begin{equation}
    \left(  \mathbf{E}\big[|F|^{p}\big]\right) ^{1/p}\leqslant c_{p,q}\left(  \mathbf{E}\big[|F|^{2}%
    \big]\right) ^{1/2}\ \mbox{ for any }p\geq 2.  \label{hypercontractivity}
    \end{equation}%
    It should be noted that the constants $c_{p,q}$ above are known with some
    precision when $F$ is a single chaos term: indeed, by Corollary 2.8.14 in
    \cite{NPbook}, $c_{p,q}=\left( p-1\right) ^{q/2}$.

\subsection{Malliavin derivative}
Let $(B_{t})_{t \in [0, T]} $ be a Wiener process  and let $\mathcal{S}$ be the class of smooth functionals of the form
\begin{equation}
\label{smooth}
F= f( B_{t_{1}},.., B_{t_{n}}), \hskip0.5cm t_{1},.., t_{n} \in [0, T],
\end{equation}
with $f \in C^{\infty} (\mathbb{R} ^{n})$ with at most polynomial growth (for $f$ and its derivatives).  For the random variable (\ref{smooth}) we define its Malliavin derivative with respect to $B$ by 
\begin{equation*}
D_{t}F= \sum_{i=1}^{n} \frac{\partial f}{\partial x_{i}} (  B_{t_{1}},.., B_{t_{n}})1_{[0, t_{i}] } (t), t \in (0, T].
\end{equation*}
The operator $D$ is an unbounded closable operator and it can be extended to the closure of $\mathcal{S} $ with respect to the norm 
\begin{equation*}
\Vert F\Vert _{k, p} ^{p} = \mathbf{E}  \vert F \vert ^{p} +\sum_{j=1} ^{k} \mathbf{E}  \Vert D ^{(j)} F \Vert ^{p} _{ L ^{2} ([0, T] ^{j})}, \hskip0.3cm F\in\mathcal{S}, p\geq 2, k\geq 1,
\end{equation*}
denoted by $\mathbb D^{k,\,p}$.

We denote by $D^{(j)}$ the $j$th iterated Malliavin derivative. The Skorohod integral integral, denoted by $\delta$, is the adjoint operator of $D$. Its domain is 
\begin{equation*}
Dom (\delta)= \left\{ u \in L ^{2} \left([0, T]\times \Omega\right), \mathbf{E} \left|  \int_{0} ^{T} u_{s} D_{s}F ds \right| \leq C \Vert F\Vert _{2} \right\}  
\end{equation*}
and we have the duality relationship
\begin{equation*}
\mathbf{E}  F \delta (u)= \mathbf{E}  \int_{0} ^{T} D_{s}F u_{s} ds , \hskip0.3cm F \in \mathcal{S}, u\in Dom (\delta).
\end{equation*}
We set $\mathbb{L} ^{k,p} = L^{p} ([0, T]; \mathbb{D} ^{k, p}), k\geq 1, p\geq 2$. This set is a subset of $\Dom (\delta)$.

\end{document}